\documentclass[12pt]{article}
\usepackage{amsfonts,amssymb}
\newcommand{\ZZ}{{\mathbb Z}}

\newcommand{\sfrac}[2]{{\textstyle{\frac{#1}{#2}}}}
\newcommand{\shalf}{{\textstyle{\frac{1}{2}}}}
\newcommand{\half}{{\frac{1}{2}}}
\newcommand{\finproof}{{\hfill \rule{5pt}{5pt}}}

\newtheorem{thm}{Theorem}
\newtheorem{lemm}{Lemma}

\begin{document}
\newpage
\pagestyle{empty}
\setcounter{page}{0}
\vfill
\begin{center}

{\Large {\bf Central extensions of classical 

\vspace{3mm}

and quantum $q$-Virasoro algebras}}

\vspace{10mm}

{\large J. Avan}

\vspace{4mm}

{\em LPTHE, CNRS-URA 280, Universit{\'e}s Paris VI/VII, France}

\vspace{7mm}

{\large L. Frappat, M. Rossi, P. Sorba}

\vspace{4mm}

{\em Laboratoire d'Annecy-le-Vieux de Physique Th{\'e}orique LAPTH, CNRS-URA 
1436} 

{\em LAPP, BP 110, F-74941 Annecy-le-Vieux Cedex, France}

\end{center}

\vfill
\vfill

\begin{abstract}
We investigate the central extensions of the $q$-deformed 
(classical and quantum) Virasoro algebras constructed from the elliptic 
quantum algebra ${\cal A}_{q,p}(sl(N))_{c}$.  After establishing the 
expressions of the cocycle conditions, we solve them, both in the 
classical and in the quantum case (for $sl(2)$).  We find that the 
consistent central extensions are much more general that those found 
previously in the literature.
\end{abstract}

\vfill
\vfill

\rightline{LAPTH-680/98}
\rightline{PAR-LPTHE 98-28}
\rightline{math.QA/9806065}
\rightline{May 1998}

\newpage
\pagestyle{plain}

%%%%%%%%%%%%%%%%%%%%%%%%%%%%%%%%%%%%%%%%%%%%%%%%%%%%%%%%%%%%%%%%%%%%%%%
\section{Introduction}

In three previous papers \cite{AFRS1,AFRS2,AFRS3} we have given a 
construction of $q$-deformed classical and quantum Virasoro and ${\cal 
W}_{N}$ algebras.  Such algebraic structures have been previously 
defined using explicit representations by $q$-deformed bosonic 
operators, associated to the collective variable representation of the 
relativistic Ruijsenaars-Schneider model \cite{SKAO,AKOS} or 
alternatively to the $q$-deformed version of the Miura transformation 
for the quantum group ${\cal U}_{q}(sl(N))_{c}$ \cite{FR,FF,BowPil}.  
These explicit constructions yielded centrally extended $q$-deformed 
algebras due to particular representation-dependent relations arising 
between their generating operators.  By contrast, our construction 
relied on abstract algebraic relations stemming from the newly proposed 
elliptic quantum algebra ${\cal A}_{q,p}(sl(N))_{c}$ \cite{JKM,JKOS}.  
The essential feature of our construction was the understanding of the 
crucial role played by the supplementary parameter $p$ (elliptic nome).  
It provided us with the intermediate step in the process of constructing 
a Poisson bracket on an extended center of ${\cal A}_{q,p}(sl(N))_{c}$, 
realizing  immediate quantizations of this Poisson structure inside the 
original algebraic object ${\cal A}_{q,p}(sl(N))_{c}$.  Elliptic 
algebras were thus shown to be a natural setting for generic 
construction of quantum $q$-${\cal W}_{N}$ algebras.

This procedure however does not naturally give rise to centrally 
extended algebras and the problem must therefore be considered 
separately.  Starting from the abstract, non-centrally extended 
quadratic structure, we shall define the general cocycle condition for 
central extensions and systematically look for its solutions.

As a first step we wish to present here a number of results concerning 
the easiest case of $q$-deformed Virasoro algebra.  We shall first of 
all describe the classical cocycle equations and give explicit sets of 
solutions to them.  At this time, we only have explicit solutions for 
the $k=0$ sector, equivalent to the original classical construction 
\cite{FR}.  They are much more general than the central extension given 
in \cite{FR} and we give a scheme of realization for them.

We then describe the quantum cocycle equations and give the general 
solution for the sector $k=0$.  This sector corresponds to an exchange 
function which is the square of the original exchange function in 
\cite{AKOS,FF}.  Interestingly enough, the cocycle condition is 
expressed as a residue formula for a particular meromorphic function!

We finally make some comments and conclusive remarks.

%%%%%%%%%%%%%%%%%%%%%%%%%%%%%%%%%%%%%%%%%%%%%%%%%%%%%%%%%%%%%%%%%%%%%%%
\section{Classical case}
\setcounter{equation}{0}

\subsection{General form of central extensions}

We consider a general quadratic Poisson algebra
\begin{equation}
\{ s_n , s_m \}=\sum _{l\in \ZZ} c_l \, s_{n-2l} \, s_{m+2l} \quad
\label{eq11}
\end{equation}
where $n,m \in 2\ZZ$ and because of antisymmetry:
\begin{equation}
c_l=-c_{-l} \,.
\label{eq12}
\end{equation}
Associativity (Jacobi) property is implied by antisymmetry.  Indeed, 
imposing:
\begin{equation}
\{ s_n , \{ s_m,s_r \} \} +\{ s_m , \{ s_r,s_n \} \} +\{ s_r, \{ s_n, 
s_m \} \} =0 \,,
\label{eq13}
\end{equation}
one obtains the conditions (where $\forall \, n,m,r \in 2\ZZ$):
\begin{eqnarray}
&&\sum _{l,j \in \ZZ} c_lc_j \left (s_{n-2j}s_{m-2l+2j}s_{r+2l} + 
s_{m-2l}s_{n-2j}s_{r+2l+2j} + s_{m-2j}s_{r-2l+2j}s_{n+2l} \right.
\nonumber \\
&& \left. + s_{r-2l}s_{m-2j}s_{n+2l+2j} + s_{r-2j}s_{n-2l+2j}s_{m+2l} + 
s_{n-2l}s_{r-2j}s_{m+2l+2j} \right) = 0 \,,
\label{eq13b}
\end{eqnarray}
which can be rewritten as:
\begin{equation}
\sum _{l,j \in \ZZ} s_{n-2j }s_{m-2l+2j} s_{r+2l} (c_lc_j + c_{l-j}c_j + 
c_{-j}c_{l-j} + c_{-l}c_{l-j} + c_{j-l}c_{-l} + c_jc_{-l}) = 0  \,.
\label{eq14}
\end{equation}
Equations (\ref {eq14}) are satisfied because of (\ref {eq12}).

\medskip

We now consider the general centrally extended version of algebra 
(\ref{eq11}):
\begin{equation}
\{ s_n, s_m \}=\sum _{l\in \ZZ}c_ls_{n-2l} s_{m+2l} + h_{n,m}
\label{eq15}
\end{equation}
where $h_{n,m}$ Poisson-commutes with $s_r$ and satisfies the 
antisymmetry property:
\begin{equation}
h_{n,m} = -h_{m,n} \,.
\label{eq16}
\end{equation}
The cocycle $h_{n,m}$ can not be a trivial one (a coboundary). Indeed 
coboundary terms are generated by the redefinitions of the generators 
$s_{n}$.  Since the bracket is quadratic, a constant (non-dynamical) 
extension is generated specifically by shifting $s_{n}$ by a constant 
$\delta_{n}$.  However, this would also simultaneously produce a linear 
term in (\ref{eq15}), and no further redefinition of $s_{n}$ by linear 
or higher terms may cancel such a linear term.
 
Condition (\ref{eq13}) now reads for all $m,n,r \in 2\ZZ$ (terms 
involving products of two $c_l$'s give no contributions because of 
(\ref{eq14})):
\begin{equation}
\sum_{l \in 2\ZZ} \left[ c_{\frac {l-r}{2}}(h_{n,m-l+r}-h_{m,n-l+r}) + 
c_{\frac {l-m}{2}}(h_{r,n-l+m}-h_{n,r-l+m}) + c_{\frac 
{l-n}{2}}(h_{m,r-l+n}-h_{r,m-l+n})\right] s_l =0 \,.
\label{eq17}
\end{equation}
Condition (\ref{eq17}) must be satisfied for all $s_l$, leading to the 
final cocycle equation for $h_{n,m}$ ($m,n,r,l \in 2\ZZ$):
\begin{equation}
c_{\frac {l-r}{2}}(h_{n,m-l+r}-h_{m,n-l+r}) + 
c_{\frac {l-m}{2}}(h_{r,n-l+m}-h_{n,r-l+m}) + 
c_{\frac {l-n}{2}}(h_{m,r-l+n}-h_{r,m-l+n}) = 0 \,. 
\label{eq18}
\end{equation}
To solve it, we note that in the particular case:
\begin{equation}
m=l=0 
\label{eq19}
\end{equation}
this equation implies:
\begin{equation}
h_{n,r}=h_{0,n+r}\frac {c_{\frac {r}{2}}-c_{\frac {n}{2}}}
{c_{\frac {r}{2}}+c_{\frac {n}{2}}} \,.
\label{eq110}
\end{equation}
Inserting (\ref {eq110}) in (\ref {eq18}), one obtains, for all $m,n,r,l 
\in 2\ZZ$:
\begin{eqnarray}
&& h_{0,n+m+r-l} \left [ c_{\frac {l-r}{2}} \displaystyle \left(
    \frac{c_{\frac{m-l+r}{2}} - c_{\frac{n}{2}}} 
         {c_{\frac{m-l+r}{2}} + c_{\frac{n}{2}}} 
  - \frac{c_{\frac{n-l+r}{2}} - c_{\frac{m}{2}}}
         {c_{\frac{n-l+r}{2}} + c_{\frac{m}{2}}} \right) +
    c_{\frac{l-m}{2}} \displaystyle \left( 
    \frac{c_{\frac{n-l+m}{2}} - c_{\frac{r}{2}}}
         {c_{\frac{n-l+m}{2}} + c_{\frac{r}{2}}} 
  - \frac{c_{\frac{r-l+m}{2}} - c_{\frac{n}{2}}}
         {c_{\frac{r-l+m}{2}} + c_{\frac{n}{2}}} \right) + 
    \right. \nonumber \\ 
&&  + \left . c_{\frac{l-n}{2}} \displaystyle \left( 
    \frac{c_{\frac{r-l+n}{2}} - c_{\frac{m}{2}}}
         {c_{\frac{r-l+n}{2}} + c_{\frac{m}{2}}} 
  - \frac{c_{\frac{m-l+n}{2}} - c_{\frac{r}{2}}}
         {c_{\frac{m-l+n}{2}} + c_{\frac{r}{2}}} \right) \right] = 0 \,. 
\label{eq111}
\end{eqnarray}
Equations (\ref{eq110}, \ref{eq111}) are thus equivalent to 
(\ref{eq18}).  Let us concentrate on (\ref{eq111}); two 
possibilities arise : \\
$i)$ the $c$-dependent factor is different from zero for every value of 
$m+n+r-l$; in this case $h_{0,M}=0$ for all $M$; using equation 
(\ref{eq110}) we obtain $h_{n,m}=0$ for all $n,m$, then the algebra 
(\ref{eq11}) cannot be centrally extended.  \\
$ii)$ the $c$-dependent factor is equal to zero for some values of
$m+n+r-l$; in this case $h_{0,M}$ is left arbitrary for these values
of $M$; using (\ref{eq110}) we then find 
the general non-zero central extensions of (\ref{eq11}). 

\medskip

Turning to our particular case we consider now the family of quadratic 
Poisson algebras indexed by a non-negative integer $k$:
\begin{equation}
\{ s_n , s_m \}_{(k)} = \sum_{l\in \ZZ} c_l^{(k)} s_{n-2l} s_{m+2l} \,, 
\quad (n,m \in 2\ZZ) \,,
\label{eq112}
\end{equation}
with
\begin{equation}
c_l^{(k)} = (-1)^{k+1} 2 \ln q \, \frac{q^{(2k+1)l}-q^{-(2k+1)l}}{q^l+q^{-l}} \,.
\label{eq112bis}
\end{equation}
It is easy to check that in the case $k=0$ the $c$-dependent factor in 
(\ref {eq111}) vanishes for every $m,n,r,l$; then $h_{0,m}$ can take arbitrary 
values for every $m$ and equation (\ref {eq110}) gives:
\begin{equation}
h_{n,m}=h_{0,n+m} \frac{q^{\frac{m-n}{2}}-q^{\frac{n-m}{2}}}
{q^{\frac {m+n}{2}}-q^{-\frac{m+n}{2}}} \,. 
\label{eq113}
\end{equation}
Grouping all the $n+m$ depending factors, one concludes that algebra
(\ref {eq112}) in the case $k=0$ has central extensions of the form:  
\begin{equation}
h_{n,m} = \xi_{n+m} \left(q^{\frac{1}{2}(m-n)} - q^{\frac{1}{2}(n-m)} \right) \,,
\label{eq114}
\end{equation}
where $\xi_n$ is a completely arbitrary function. 

Remark that the particular choice $\xi_n= -2 \ln q \, \delta _{n,0}$ 
gives the central extension found by Frenkel and Reshetikhin \cite{FR} 
for the Poisson bracket (\ref{eq112}) at $k=0$.  However our 
calculations show that this Poisson bracket allows for more general 
non-local central extensions of the form (\ref{eq114}).

\medskip

The situation is not fully clarified for $k \not = 0$.  Equation 
(\ref{eq111}) implies in fact by numerical calculations $h_{0,0} = 
h_{0,1} = h_{0,2} = 0$, but we have not yet concluded for or against the 
existence of higher label $h_{0,n}$ central extensions ($n > 2$).

\subsection{Explicit construction of central extensions}

We have showed that algebra (\ref{eq112}) in the case $k=0$ 
can be centrally extended in the following way:
\begin{equation}
\{ s_n, s_m \} = -2 \ln q \, \sum _{l\in \ZZ}\frac
{q^l-q^{-l}}{q^l+q^{-l}}s_{n-2l}s_{m+2l} + \xi_{m+n} \left( 
q^{\frac{m-n}{2}}-q^{\frac {n-m}{2}} \right) \quad , \quad m,n \in 2\ZZ \, .
\label{eq115}
\end{equation}
We now give an explicit realization of (\ref{eq115}) in terms of 
elements belonging to the center of ${\cal A}_{q,p}(\widehat{sl}(2)_c)$ 
at $c=-2$.

\medskip

{From} the result of \cite{AFRS1} we know that the operators:
\begin{equation}
t(z)={\mbox {Tr}}\left( L^+(q^{\frac {c}{2}}z)L^-(z)^{-1}\right) \,
\label{eq116}
\end{equation}
satisfy the Poisson algebra at $c = -2$:
\begin{equation}
\Big\{ t(z),t(w) \Big\} = {\cal Y}(z/w) \, t(z) \, t(w) \, ,
\label{eq117}
\end{equation}
where
\begin{eqnarray}
{\cal Y}(x) &=& -(2\ln q) \left[ \sum_{n \ge 0} ~ \left( 
\frac{2x^2q^{4n+2}}{1-x^2q^{4n+2}} - \frac{2x^{-2}q^{4n+2}}{1-x^{-2}q^{4n+2}} 
\right) + \right. \nonumber \\
&& \left. +\sum_{n>0} ~ \left( 
- \frac{2x^2q^{4n}}{1-x^2q^{4n}} + \frac{2x^{-2}q^{4n}}{1-x^{-2}q^{4n}}
\right) - \frac{x^2}{1-x^2} + \frac{x^{-2}}{1-x^{-2}} \right] \,\,.
\label{eq118}
\end{eqnarray}
Because of the property:
\begin{equation}
{\cal Y}(x) = -{\cal Y}(xq) \,,
\label{eq119}
\end{equation}
the field
\begin{equation}
h(z)=t(zq)
\label{eq120}
\end{equation}
satisfies the Poisson brackets:
\begin{eqnarray}
\Big\{ h(z),t(w)\Big\}&=&-{\cal Y}(z/w) \, h(z) \, t(w)  \nonumber \\
\Big\{ h(z),h(w)\Big\}&=&{\cal Y}(z/w) \, h(z) \, h(w) \,.
\label{eq121}
\end{eqnarray}
Writing equations (\ref {eq117}),(\ref {eq121}) in terms of the modes:
\begin{equation}
t_n = \oint_C \frac{dz}{2\pi iz} \, z^{-n} \, t(z) \,, \quad 
h_n = \oint_C \frac{dz}{2\pi iz} \, z^{-n} \, h(z) \,,
\label{eq121b}
\end{equation}
and extending as previously done the validity of (\ref{eq117}) to the level of 
analytic continuations of the structure function ${\cal Y}(z/w)$,
we obtain \cite {AFRS1} three families of Poisson brackets depending on a 
non-negative integer $k$:
\begin{eqnarray}
\{ t_n,t_m \} &=& \sum_{l\in \ZZ} c^{(k)}_l \, t_{n-2l} \, t_{m+2l}
\nonumber \\
\{ h_n,t_m  \} &=& -\sum _{l\in \ZZ} c^{(k)}_l \, h_{n-2l} \, t_{m+2l}
\label{eq122} \\
\{ h_n,h_m  \} &=& \sum _{l\in \ZZ} c^{(k)}_l \, h_{n-2l} \, h_{m+2l} \,,
\nonumber
\end{eqnarray}
where $c^{(k)}_l$ is given by equation (\ref {eq112bis}). We are interested
in the case $k=0$: we define the observable
\begin{equation}
s_n = q^{-\frac{n}{2}} t_n + q^{\frac{n}{2}} h_n 
\label{eq123}
\end{equation}
Using (\ref {eq122}), the Poisson bracket now reads:

\begin{eqnarray}
\{ s_n,s_m  \}&=& -2 \ln q \, \Big[ \sum _{l \in \ZZ} \frac{q^l-q^{-l}}{q^l+q^{-l}} 
\, s_{n-2l} \, s_{m+2l}-\nonumber \\
&&- \sum _{l\in \ZZ} (1-q^{-2l}) q^{\frac{n-m}{2}} 
\, h_{n-2l} \, t_{m+2l} + \sum _{l \in \ZZ} (1-q^{2l}) 
q^{\frac{m-n}{2}} \, h_{m+2l} \, t_{n-2l} \Big] \,. 
\label{eq126}
\end{eqnarray}
Remark that the analytic simplifications leading from (\ref{eq122}) to 
(\ref{eq126}) are specific of the behaviour of the 
structure coefficients $c_{l}^{(k)}$ 
at $k=0$, which gives a technical justification for
the (easily proved) existence of central extensions at 
$k=0$ and the difficulty to obtain central extensions at $k \ne 0$. \\
Performing the shifts:
\begin{equation}
l\rightarrow l-\frac {m}{2} \,, \quad l\rightarrow \frac {n}{2}-l
\, , \label{eq127}
\end{equation}
respectively in the penultimate and in the last term of (\ref{eq126}), 
we can simplify them:
\begin{equation}
\{ s_n,s_m \} = -2 \ln q \, \Big[ \sum_{l\in \ZZ} 
\frac{q^l-q^{-l}}{q^l+q^{-l}} \, s_{n-2l} \, s_{m+2l} + \left( 
q^{\frac{m-n}{2}}-q^{\frac{n-m}{2}} \right) \sum_{l\in \ZZ} h_{n+m-2l} 
\, t_{2l} \Big] \,.
\label{eq128}
\end{equation}
We now prove that the last term is central in the algebra generated by 
$t_n$ and $h_n$, that is:
\begin{eqnarray}
\{t_r, \sum_{l\in \ZZ} h_{n+m-2l} \, t_{2l} \} &=& 0 \nonumber \\
\{h_r, \sum_{l\in \ZZ} h_{n+m-2l} \, t_{2l} \} &=& 0 \,.
\label{eq129}
\end{eqnarray}
Let us consider the first Poisson bracket.  The use of (\ref{eq122}) for 
$k=0$ gives:
\begin{equation}
\{t_r, \sum_{l\in \ZZ} h_{n+m-2l} \, t_{2l} \} = 
-2 \ln q \, \sum_{l,l^\prime \in \ZZ} \left( 
\frac{q^{l^\prime}-q^{-l^\prime}}{q^{l^{\prime}}+q^{-l^\prime}} \,
h_{n+m-2l-2l^{\prime}} \, t_{r+2l^{\prime}}t_{2l} + 
\frac{q^{l^\prime}-q^{-l^\prime}}{q^{l^\prime}+q^{-l^\prime}} \, 
h_{n+m-2l}t_{r-2l^{\prime}} \, t_{2l+2l^{\prime}} \right) \,. 
\label{eq130}
\end{equation}
Making the changes of variables $l^{\prime}\rightarrow -l^{\prime}$ and 
$l\rightarrow l+l^{\prime}$ in the last term of (\ref {eq130}), one 
obtains immediately the first identity of (\ref {eq129}).  The second 
one is proved in an identical way.

\medskip

Any centrally extended algebra (\ref{eq15})-(\ref{eq114}) may thus 
be generated by the combination of the abstract generators $t_{n}$, $h_{m}$ in 
(\ref{eq122}), after setting each constant quantity $-2 \ln q \, \sum_{l 
\in \ZZ} h_{n-2l} \, t_{2l}$ to a fixed value $\xi_{n}$.  Notice that 
the Frenkel-Reshetikhin construction uses a representation of $t$ and 
$h$ such that $h(z) \equiv t(z)^{-1}$, hence indeed $\xi_{n} = -2 \ln q 
\, \delta_{n,0}$.

%%%%%%%%%%%%%%%%%%%%%%%%%%%%%%%%%%%%%%%%%%%%%%%%%%%%%%%%%%%%%%%%%%%%%%%
\section{Quantum case}
\setcounter{equation}{0}

We start from the general symmetrized exchange relation \cite{AFRS2,AFRS3}
\begin{equation}
\sum_{l \in 2\ZZ} f_{l} (t_{n-l} t_{m+l} - t_{m-l} t_{n+l}) = 0 \,.
\label{eq31}
\end{equation}
The structure constants $f_{l}$ are given by \cite{AFRS3} (hence the
denomination``symmetrized''):
\begin{equation}
f_{l} = \shalf \, (f_{l}^{(0)} + f_{l}^{(1)}) \,,
\label{eq31ter}
\end{equation}
where $f_{l}^{(0)}$ and $f_{l}^{(1)}$ are respectively the coefficients 
of the series expansion of the meromorphic function $f(z)$ inside the 
convergence 
disk $\vert z \vert < 1$ and in the convergence ring $1 < \vert z \vert 
< \vert p^{-\half}q^{-1} \vert$ with
\begin{equation}
f(z) = \frac{1}{(1-z^2)^2} \prod_{n \ge 0} 
\frac{(1-z^2pq^{4n})^2(1-z^2p^{-1}q^{4n+2})^2}
{(1-z^2pq^{4n+2})^2(1-z^2p^{-1}q^{4n+4})^2} \,.
\label{eq31bis}
\end{equation}
Equation (\ref{eq31}) can be written as:
\begin{equation}
\Big[ t_{n} , t_{m} \Big] = \sum_{l \in 2\ZZ^*} -f_{l} (t_{n-l} t_{m+l} 
- t_{m-l} t_{n+l}) \,.
\label{eq32}
\end{equation}
Hence a centrally extended version of (\ref{eq32}) must take the form
\begin{equation}
\Big[ t_{n} , t_{m} \Big] = \sum_{l \in 2\ZZ^*} -f_{l} (t_{n-l} t_{m+l} 
- t_{m-l} t_{n+l}) + h_{n,m} \,.
\label{eq33}
\end{equation}
In order to determine the solutions of (\ref{eq33}), one writes the 
Jacobi identity:
\begin{equation}
\Big[ t_{n} , \Big[ t_{m} , t_{r} \Big] \Big] + 
\Big[ t_{m} , \Big[ t_{r} , t_{n} \Big] \Big] + 
\Big[ t_{r} , \Big[ t_{n} , t_{m} \Big] \Big] = 0 \,.
\label{eq34}
\end{equation}
Inserting (\ref{eq33}) into (\ref{eq34}), one obtains a condition 
involving trilinear terms of the form $t_{m-k}h_{r,n+k}f_{k}$ and 
quartic terms of the form $t_{m-k-l}h_{r+l,n+k}f_{k}f_{l}$, summed over 
$k,l \in 2\ZZ^{*}$.  These sums can be rewritten  by summing 
over the indices $k,l \in 2\ZZ$, $f_{k}$ being replaced by $f'_{k} = 
f_{k} - \delta_{k,0}$.  Plugging the explicit expression of $f'_{k}$, 
the sum with the trilinear terms vanishes and one is left with the 
cocycle condition:
\begin{equation}
\sum_{k \in 2\ZZ} (f_{k-n} f_{s-r-k} - f_{n-k} f_{k+r-s} + f_{k-m} 
f_{s-n-k} - f_{m-k} f_{k+n-s} + f_{k-r} f_{s-m-k} - f_{r-k} f_{k+m-s} ) 
\, h_{s-k,k} = 0 \,.
\label{eq35}
\end{equation}
which has to be satisfied for all possible triplets $n,m,r \in 2\ZZ$ 
and for all values of $s \in 2\ZZ$.  \\
We will first focus on solutions of the form $h_{n,m} = \lambda_{m} 
\delta_{n+m,0}$, i.e. the subset of cocycle equations
for $s=0$.  The cocycle condition (\ref{eq35}) becomes:
\begin{equation}
\sum_{k \in 2\ZZ} ( f_{k-n} f_{-r-k} - f_{n-k} f_{k+r} + f_{k-m} 
f_{-n-k} - f_{m-k} f_{k+n} + f_{k-r} f_{-m-k} - f_{r-k} f_{k+m} ) 
\, \lambda_{k} = 0 \,.
\label{eq35bis}
\end{equation}

The resolution of the cocycle condition will be done in four steps.  In 
the first step, we show that if the solution of (\ref{eq35bis}) exists, 
it is necessary of the form $\lambda_{k} = q^{k}-q^{-k}$.  In the second 
step, we prove that the series expansions that appears in 
(\ref{eq35bis}) can then be resummed consistently.  In the third step, 
we transform the derived equation into a contour integral that allows 
to achieve the proof.  In the last step, we establish the immediate 
generalization of the form of the cocycle term for $s \neq 0$ in 
(\ref{eq35bis}).

\bigskip

\noindent \textbf{Step 1}: \\
The form of the coefficients $\lambda_{k}$ can be guessed as follows.  
Consider all possible triplets $(n,m,r)$ such that $\max(n,m,r) \le M$ 
where $M$ is a \emph{fixed} integer and write all cocycle conditions 
(\ref{eq35bis}) corresponding to these triplets. One assumes
of course that (\ref{eq35bis}) are absolutely convergent series
which can be manipulated consistently. This will then be proved
in Step 2 in order to guarantee the consistency of the overall derivation.

 It appears that the 
obtained relations differ only by finite numbers of terms.  One gets in 
this way a finite set of equations, each one having a finite number of 
terms, involving the coefficients $\lambda_{k}$ with $\vert k \vert \le 
M$.  It is then possible to solve these equations step by step.  
Numerical calculations for $M = 4,6,8$ show that the solution is 
necessarily of the form $\lambda_{k} = q^{k}-q^{-k}$ (up to an overall 
factor).  \\
Therefore, we only need to show that the cocycle condition (\ref{eq35bis}) 
with the choice $\lambda_{k} = q^{k}-q^{-k}$ for all $k \in 2\ZZ$ is 
satisfied, that is:
\begin{eqnarray}
&& \sum_{k \in 2\ZZ} (f^0_{k-n} f^0_{-r-k} - f^0_{n-k} f^0_{k+r}
+ f^0_{k-m} f^0_{-n-k} - f^0_{m-k} f^0_{k+n} 
+ f^0_{k-r} f^0_{-m-k} - f^0_{r-k} f^0_{k+m} ) \, (q^{k}-q^{-k}) + 
\nonumber \\ 
&& \sum_{k \in 2\ZZ} \alpha \Big( (m-n)(f^0_{k+r}+f^0_{-k-r}) + 
(r-m)(f^0_{k+n}+f^0_{-k-n}) + (n-r)(f^0_{k+m}+f^0_{-k-m}) \Big)
\, (q^{k}-q^{-k}) = 0 \nonumber \\
\label{eq36}
\end{eqnarray}
where $\alpha = -\sfrac{1}{4} (1-x^2)^2 f(x)\Big\vert_{x=1}$.  

\medskip

\noindent
For convenience, we denote by $S_{1}$ and $S_{2}$ respectively the 
first and second sums in (\ref{eq36}). Remark that $S_{1}$ is in fact a 
finite sum since $f^0_{l} = 0$ for $l < 0$.

\bigskip

\noindent \textbf{Step 2}: \\

\begin{lemm}\label{lemm1}
The sum $S_{2}$ is a convergent series expansion, the sum of which is 
equal to
\begin{equation}
S_{2} = -\alpha \Big(f(q) + f(q^{-1})\Big) \Big((m-n)(q^{r}-q^{-r}) + 
(n-r)(q^{m}-q^{-m}) + (r-m)(q^{n}-q^{-n})\Big) \,.
\label{eq37}
\end{equation}
\end{lemm}
\textbf{Proof}: \\
One has
\begin{eqnarray}
S_{2} &=& \sum_{cyclic(m,n,r)} \alpha (m-n) \sum_{k \in 2\ZZ} q^{k} \, 
f^0_{k+r} - \sum_{cyclic(m,n,r)} \alpha (m-n) \sum_{k \in 2\ZZ} q^{-k} 
\, f^0_{-k-r} \nonumber \\
&& + \sum_{cyclic(m,n,r)} \alpha (m-n) \sum_{k \in 2\ZZ} 
(q^{k} \, f^0_{-k-r} - q^{-k} \, f^0_{k+r}) \,.
\end{eqnarray}
The first two sums are convergent and their sum is equal to 
$\displaystyle \sum_{cyclic(m,n,r)} -\alpha (m-n)(q^{r}-q^{-r})f(q)$.  
To study the last sum, let $f^1_{l}$ be the coefficients of the series 
expansion of $f(x)$ in the convergence ring $1 < \vert x \vert < \vert 
p^{-\half}q^{-1} \vert$.  One has $f^1_{l} = f^0_{l} + 2\alpha l + 
\beta$ where $\beta$ is a constant, $\alpha$ is defined above and $l \in 
\ZZ$.  The last sum can then be rewritten as
\begin{eqnarray}
&& \sum_{cyclic(m,n,r)} \alpha (m-n) \sum_{k \in 2\ZZ} \Big(q^{k} \, 
(f^1_{-k-r} + 2\alpha (k+r) - \beta) - q^{-k} \, (f^1_{k+r} - 2\alpha 
(k+r) - \beta)\Big) \nonumber \\
&& \qquad = \sum_{cyclic(m,n,r)} \alpha (m-n) \sum_{k \in 2\ZZ} (q^{k} \, 
f^1_{-k-r} - q^{-k} \, f^1_{k+r} ) \nonumber \\
&& \qquad \quad + \sum_{cyclic(m,n,r)} \alpha (m-n) \sum_{k \in 2\ZZ} 
\Big( 2\alpha (k+r)(q^{k} + q^{-k}) - \beta(q^{k} - q^{-k}) \Big) \,.
\end{eqnarray}
The first sum is convergent and is equal to $\displaystyle 
\sum_{cyclic(m,n,r)} -\alpha (m-n)(q^{r}-q^{-r})f(q^{-1})$, while the 
second sum is identically zero term by term  due to cyclicity over the indices 
$m,n,r$ since $\displaystyle \sum_{cyclic(m,n,r)} (m-n) = 0$ and 
$\displaystyle \sum_{cyclic(m,n,r)} (m-n)r = 0$. 

It is understood at every
step of the derivation that the cyclization must be done before the summation 
over $k \in 2\ZZ$. 

This achieves the 
proof of Lemma \ref{lemm1}. 
\finproof

\medskip

\noindent
Note that the definition of $\alpha$ (see above) and the expression of 
$f(z)$ (see (\ref{eq31bis})) implies that
\begin{equation}
-\alpha \Big(f(q) + f(q^{-1})\Big) = \shalf \, 
\frac{(p-1)^{2}(p-q^{2})^{2}}{p^{2}(1-q^{2})^{2}} \equiv \shalf \, 
{\cal E} \,.
\label{eq37bis}
\end{equation}

\bigskip

\noindent \textbf{Step 3}: \\
\begin{lemm}
The cocycle condition (\ref{eq35}) is equivalent to the following 
contour integral, the contour enclosing all the poles (including the 
pole at infinity) of the integrand (that is $0$, $\pm q$, $\pm 1$, 
$\pm q^{-1}$, $\infty$):
\begin{eqnarray}
&\displaystyle \oint \frac{dz}{2i \pi z} \bigg[ & 
q^{-r} z^{r+n} f(z) f(zq^{-1}) - q^{r} z^{r+n} f(z) f(zq)
+ q^{r} z^{-r-n} f(z) f(zq^{-1}) \nonumber \\
&& - q^{-r} z^{-r-n} f(z) f(zq) + q^{-m} z^{m+r} f(z) f(zq^{-1}) 
- q^{m} z^{m+r} f(z) f(zq) \bigg. \nonumber \\
&& + q^{m} z^{-m-r} f(z) f(zq^{-1}) - q^{-m} z^{-m-r} f(z) f(zq)
+ q^{-n} z^{n+m} f(z) f(zq^{-1}) \bigg. \nonumber \\
&& - q^{n} z^{n+m} f(z) f(zq) + q^{n} z^{-n-m} f(z) f(zq^{-1}) 
- q^{-n} z^{-n-m} f(z) f(zq) \ \ \bigg] = 0 \,. \nonumber \\
\label{eq310}
\end{eqnarray}
\end{lemm}
\textbf{Proof:} Let us introduce the function $g_{\pm}(z)$ defined by:
\begin{equation}
g_{\pm}(z) \equiv \frac{(1-p^{\pm 1}z)^2(1-p^{\mp 1}q^{\pm 2}z)^2} 
{(1-z)^2(1-q^{\pm 2}z)^2} \,,
\label{eq311}
\end{equation}
which satisfies
\begin{equation}
f(z) f(zq^{\pm 1}) = g_{\pm}(z^2) \,.
\label{eq312}
\end{equation}
By virtue of (\ref{eq311})-(\ref{eq312}), the integrand in (\ref{eq310}) 
has double poles at $z = 0$, $z = \pm 1$, $z = \pm q^{\pm 1}$, $z = 
\infty$.  One has
\begin{eqnarray}
\oint_{C_{0}} \frac{dz}{2i \pi z} z^{r+n} (q^{-r} f(z) f(zq^{-1}) - 
q^{r} f(z) f(zq)) &=& \oint \frac{dz}{2i \pi z} \sum_{l \ge 
0} \sum_{l' \ge 0} f^0_{l} f^0_{l'} (q^{-l'-r} - q^{l'+r}) \, z^{r+n+l+l'} 
\nonumber \\
&=& \sum_{l \ge 0} \sum_{l' \ge 0} f^0_{l} f^0_{l'} 
(q^{-l'-r} - q^{l'+r}) \, \delta_{n+r+l+l',0} \nonumber \\
&=& \sum_{l \ge 0} f^0_{l} f^0_{-l-n-r} (q^{n+l} - q^{-n-l}) \nonumber \\ 
&=& \sum_{k \in 2\ZZ} f^0_{k-n} f^0_{-r-k} (q^{k}-q^{-k}) \,,
\label{eq313}
\end{eqnarray}
where the contour $C_{0}$ is a small circle around the origin such that 
$f(z)$ be analytic inside the contour $C_{0}$.  Summing all the terms of 
the form (\ref{eq313}) occuring in (\ref{eq310}), one sees that the 
contribution of the pole at $z = 0$ in (\ref{eq310}) is exactly equal to 
the first sum of (\ref{eq36}).  \\
Moreover, the integral (\ref{eq310}) being invariant under the change $z 
\rightarrow z^{-1}$, the contribution of the pole at infinity is equal 
to the contribution to the pole at the origin.

\medskip

\noindent
It remains to compute the contributions of the poles at $z = \pm 1$, and 
$z = \pm q^{\pm 1}$.  The calculation is greatly simplified if one 
performs the change of variables $w = z^{2}$.  One obtains:
\begin{eqnarray}
&\displaystyle \oint_{C_{z_{0}} \cup C_{-z_{0}}} \frac{dz}{2i \pi z} \bigg[ & 
q^{-r} z^{r+n} f(z) f(zq^{-1}) - q^{r} z^{r+n} f(z) f(zq)
+ q^{r} z^{-r-n} f(z) f(zq^{-1}) \nonumber \\
&& - q^{-r} z^{-r-n} f(z) f(zq) + q^{-m} z^{m+r} f(z) f(zq^{-1}) 
- q^{m} z^{m+r} f(z) f(zq) \bigg. \nonumber \\
&& + q^{m} z^{-m-r} f(z) f(zq^{-1}) - q^{-m} z^{-m-r} f(z) f(zq)
+ q^{-n} z^{n+m} f(z) f(zq^{-1}) \bigg. \nonumber \\
&& - q^{n} z^{n+m} f(z) f(zq) + q^{n} z^{-n-m} f(z) f(zq^{-1}) 
- q^{-n} z^{-n-m} f(z) f(zq) \ \ \bigg] \nonumber \\
&= \displaystyle \oint_{C_{w_{0}}} \frac{dw}{2i \pi w} \bigg[ & 
g_{-}(w) \Big( q^{-r} w^{\half(r+n)} + q^{r} w^{\half(-r-n)} + q^{-m} 
w^{\half(m+r)} + q^{m} w^{\half(-m-r)} \nonumber \\
&& + q^{-n} w^{\half(n+m)} + q^{n} w^{\half(-n-m)} \Big) + g_{+}(w) \Big(
- q^{r} w^{\half(r+n)} - q^{-r} w^{\half(-r-n)} \nonumber \\
&& - q^{m} w^{\half(m+r)} - q^{-m} w^{\half(-m-r)} - q^{n} 
w^{\half(n+m)} - q^{-n} w^{\half(-n-m)} \Big) \bigg] \,,
\label{eq314}
\end{eqnarray}
where the contours $C_{\pm z_{0}}$ are small circles around the points 
$z = \pm z_{0}$ and the contour $C_{w_{0}}$ is a small circle around the 
point $w = w_{0} = z_{0}^{2}$.  Notice that the factor $\shalf$ coming 
from the integration measure $dz/2i \pi z \rightarrow \shalf \, dw/2i 
\pi w$ is compensated by the fact that when the variable $z$ winds 
around $z_{0}$ \emph{and} $-z_{0}$ once, the variable $w$ winds around 
$w_{0}$ twice.

\medskip

\noindent 
The contribution of the poles at $w = q^{2}$ comes from the relevant
factors of $g_{-}(w)$.  It is given by:
\begin{eqnarray}
&& {\cal E} \Big[ \Big( \frac{r+n}{2} -1 \Big) q^{n} + 
\Big( \frac{m+r}{2} -1 \Big) q^{r} + 
\Big( \frac{n+m}{2} -1 \Big) q^{m} - 
\Big( \frac{m+r}{2} +1 \Big) q^{-r} \nonumber \\
&& - \Big( \frac{n+m}{2} +1 \Big) q^{-m} - 
\Big( \frac{r+n}{2} +1 \Big) q^{-n} \Big] + 
{\cal E}' \,\Big(q^{n} + q^{m} + q^{r} + q^{-n} + q^{-m} + q^{-r} \Big) \,,
\end{eqnarray}
where ${\cal E}$ is given by (\ref{eq37bis}) and ${\cal E}' = -2{\cal E} 
\displaystyle \Big( \frac{1}{p^{-1}-1} + \frac{1}{q^{-2}p-1} - 
\frac{1}{q^{-2}-1} \Big)$.  \\
In the same way, the contribution of the poles at $w = q^{-2}$ comes 
from the relevant factors of $g_{+}(w)$.  It is given by:
\begin{eqnarray}
&& {\cal E} \Big[ \Big( \frac{m+r}{2} +1 \Big) q^{r} + 
\Big( \frac{n+m}{2} +1 \Big) q^{m} + 
\Big( \frac{r+n}{2} +1 \Big) q^{n} - 
\Big( \frac{m+r}{2} -1 \Big) q^{-r} \nonumber \\
&& - \Big( \frac{n+m}{2} -1 \Big) q^{-m} - 
\Big( \frac{r+n}{2} -1 \Big) q^{-n} \Big] + 
{\cal E}'' \,\Big(q^{n} + q^{m} + q^{r} + q^{-n} + q^{-m} + q^{-r} \Big) \,,
\end{eqnarray}
where ${\cal E}'' = 2{\cal E} \displaystyle \Big( \frac{1}{p-1} + 
\frac{1}{q^{2}p^{-1}-1} - \frac{1}{q^{2}-1} \Big)$.  \\
Finally, the contribution of the poles at $w = 1$ comes from the 
relevant factors of both $g_{-}(w)$ and $g_{+}(w)$.  It is given by:
\begin{eqnarray}
&& {\cal E} \Big[ \Big( \frac{r+n}{2} -1 \Big) q^{-r} + 
\Big( \frac{m+r}{2} -1 \Big) q^{-m} + 
\Big( \frac{n+m}{2} -1 \Big) q^{-n} - 
\Big( \frac{r+n}{2} +1 \Big) q^{r} \nonumber \\
&& - \Big( \frac{m+r}{2} +1 \Big) q^{m} - 
\Big( \frac{n+m}{2} +1 \Big) q^{n} \Big] - 
{\cal E}'' \,\Big(q^{n} + q^{m} + q^{r} + q^{-n} + q^{-m} + q^{-r} \Big) 
\nonumber \\
&& + {\cal E} \Big[ - \Big( \frac{r+n}{2} -1 \Big) q^{r} - 
\Big( \frac{r+m}{2} -1 \Big) q^{m} - 
\Big( \frac{n+m}{2} -1 \Big) q^{n} + 
\Big( \frac{r+n}{2} +1 \Big) q^{-r} \nonumber \\
&& + \Big( \frac{m+r}{2} +1 \Big) q^{-m} + 
\Big( \frac{n+m}{2} +1 \Big) q^{-n} \Big] - 
{\cal E}' \,\Big(q^{n} + q^{m} + q^{r} + q^{-n} + q^{-m} + q^{-r} \Big) \,.
\end{eqnarray}
Summing up all contributions, one is left with
\begin{equation}
{\cal E} \Big[ (m-n)(q^{r}-q^{-r}) + 
(n-r)(q^{m}-q^{-m}) + (r-m)(q^{n}-q^{-n}) \Big] \,,
\end{equation}
which is exactly the expression given in eq. (\ref{eq37}) up to a factor 
2.

\medskip

\noindent 
Finally, the contour integral (\ref{eq310}) is equal to twice the 
cocycle condition (\ref{eq36}).  The function to be integrated in 
(\ref{eq310}) being meromorphic and the contour of integration 
encircling all the poles of the function (including the pole at 
infinity), one concludes that the contour integral (\ref{eq310}) is 
equal to zero.  It follows that the cocycle condition (\ref{eq36}) is 
identically satisfied.  
\finproof

\bigskip

\noindent \textbf{Step 4}: \\
We want now to relax the constraint $h_{n,m} = \lambda_{m} 
\delta_{n+m,0}$.  Let us return to the general cocycle condition 
(\ref{eq35}) which has to be satisfied for any fixed value $s \in 2\ZZ$.  
\\
$\bullet$ case 1: $s \in 4\ZZ$ \\
We perform the change of variable $m' = m - \shalf s$, $n' = n - \shalf 
s$, $r' = r - \shalf s$ and $k' = k - \shalf s$. (\ref{eq35}) is then
written as:
\begin{eqnarray}
&& \sum_{k' \in 2\ZZ} (f_{k'-n'} f_{-r'-k'} - f_{n'-k'} f_{k'+r'} + 
f_{k'-m'} f_{-n'-k'} - f_{m'-k'} f_{k'+n'+} \nonumber \\ 
&& \hspace{20mm} + f_{k'-r'} f_{-m'-k'} - f_{r'-k'} f_{k'+m'} ) 
\, h_{\half s-k',\half s+k'} = 0 \,.
\label{eq40}
\end{eqnarray}
>From the previous results, this condition is satisfied, for each value 
$s \in 4\ZZ$, by
\begin{equation}
h_{\half s-k',\half s+k'} = \xi_{s} \, (q^{k'} - q^{-k'})
\label{eq41}
\end{equation}
$\bullet$ case 2: $s \in 4\ZZ+2$ \\
The previous shift now leads to the condition
\begin{eqnarray}
&& \sum_{k' \in 2\ZZ+1} (f_{k'-n'} f_{-r'-k'} - f_{n'-k'} f_{k'+r'} + 
f_{k'-m'} f_{-n'-k'} - f_{m'-k'} f_{k'+n'} \nonumber \\
&& \hspace{20mm} + f_{k'-r'} f_{-m'-k'} - f_{r'-k'} f_{k'+m'} ) 
\, h_{\half s-k',\half s+k'} = 0 \,,
\label{eq42}
\end{eqnarray}
where $m',n', r', \shalf s \in 2\ZZ+1$.  \\
As before, the condition (\ref{eq42}) is satisfied by $h_{\half 
s-k',\half s+k'} = \xi_{s} \, (q^{k'} - q^{-k'})$ for each value $s \in 
2\ZZ+1$.  Indeed, the proof runs along the same lines as for (\ref{eq36}), 
in particular the contour integral formulation also holds.

\bigskip

\noindent
The results above can then be summarized in the following theorem:
\begin{thm}
The $q$-deformed Virasoro algebra ($m,n \in 2\ZZ$)
\[
\Big[ t_{n} , t_{m} \Big] = \sum_{l \in 2\ZZ^*} -f_{l} (t_{n-l} t_{m+l} 
- t_{m-l} t_{n+l})
\]
where the structure constants $f_{l}$ are given by 
(\ref{eq31ter})-(\ref{eq31bis}), admits a central extension of the form
\[
h_{n,m} = \xi_{n+m} \, (q^{\frac{1}{2}(n-m)} - q^{\frac{1}{2}(m-n)}) 
\]
where $\xi_{n}$ is a completely arbitrary function.
\end{thm}

We must remark that the cocycle condition (\ref{eq35}) cannot be 
fulfilled if one chooses $f_{l} = f_{l}^{0}$ instead of the 
symmetrized coefficient $\shalf(f_{l}^{0} + f_{l}^{1})$: the contour 
integral argument is not applicable in this case since the cocycle equation
contains only the sum $S_1$ which does not vanish.

In addition, the replacement of the ($k=0$) symmetrized structure 
coefficients in (\ref{eq31}) by higher label structure functions 
containing contributions from other poles of ${\cal Y}(z/w)$ leads to 
serious difficulties in the evaluation of the cocycle condition 
(\ref{eq36}) since the sum $S_{2}$ then becomes divergent owing to the
presence of quadratic and higher-power terms in the structure coefficients
$f_l$ compared to the simple linear dependance in $f_l^{(1)}$ which is
canceled by permutation arguments. 

It seems therefore that the symmetrized ($k=0$) sector for the 
structure function (\ref{eq31bis}) is the only one admitting consistent 
central extensions.

\section{Perspectives}

The results which we have obtained here, and the limitations which we 
have encountered, clearly indicate two further directions of 
investigations.  One is the study of more general central plus linear 
extensions, both classical and quantum.  The cocycle conditions are 
necessarily much more complicated, involving two sets of unknowns 
$c_{nm}^r$ (linear) and $h_{nm}$ (constants); furthermore, there now 
exist cobounbary-type extensions as follows from the comment in section 
1.

On the other hand this increase in the number of degrees of freedom relaxes
the constraints in the resolution of the cocycle equations, leading us
to hope for the existence of solutions for all structure functions.

The other direction in which spirit this investigation has been pursued 
is the search for explicit operator representations of the centrally 
extended algebras.  Based on the previous constructions 
\cite{FR,FF,BowPil} we indeed expect that such representations will 
contain central (or linear?)  extensions, motivating our current study.
  
At this time, only for the simplest structure function originally 
obtained in \cite{SKAO} were such constructions proposed, using deformed 
bosons \cite{FF,BowPil}.  Even our simplest structure function 
\cite{AFRS1}, essentially equal to the square of the structure function 
for \cite{SKAO} has yet no explicit simple operatorial representation.  
It actually seems \cite{AFRSn} that mere simple extensions by 
``normalizations'' of the vertex operator constructions in 
\cite{FF,BowPil} do not lead to realizations of this quantum algebra.  
The problem is thus fully open.

\end{document}